\documentclass{icmart}

\usepackage{psfig}
\contact[jqfan@Princeton.EDU]{Jianqing Fan\\
Department of Operation Research and Financial Engineering\\
and Bendheim Center for Finance\\
Princeton University\\
Princeton, NJ 08544}

\contact[rli@stat.psu.edu]{Runze Li\\
Department of Statistics and The Methodology Center\\
The Pennsylvania State University\\
University Park, PA 16802-2111
}

\makeatletter
\def\singlespace{\def\baselinestretch{1}\@normalsize}

\@addtoreset{equation}{section}
\renewcommand{\theequation}{\thesection.\arabic{equation}}

\renewcommand{\hat}{\widehat}
\def\singlespace{\def\baselinestretch{1}\@normalsize}

\makeatother

\newcommand{\ba}{\mbox{\bf a}}
\newcommand{\bb}{\mbox{\bf b}}

\newcommand{\bx}{\mbox{\bf x}}
\newcommand{\by}{\mbox{\bf y}}
\newcommand{\bz}{\mbox{\bf z}}

\newcommand{\bI}{\mbox{\bf I}}
\newcommand{\bW}{\mbox{\bf W}}
\newcommand{\bX}{\mbox{\bf X}}
\newcommand{\pl}{p_{\lambda}}
\newcommand{\bbeta}{\mbox{\boldmath$\beta$}}
\newcommand{\btheta}{\mbox{\boldmath$\theta$}}
\newcommand{\bxi}{\mbox{\boldmath$\xi$}}

\newcommand{\hbbeta}{\hat{\bbeta}}

\newcommand{\bla}{\mbox{\boldmath$\lambda$}}

\newcommand{\eps}{\varepsilon}

\newcommand{\bff}{\mbox{\bf f}}
\newcommand{\bA}{\mbox{\bf A}}
\newcommand{\bB}{\mbox{\bf B}}
\newcommand{\bD}{\mbox{\bf D}}
\newcommand{\bG}{\mbox{\bf G}}
\newcommand{\bL}{\mbox{\bf L}}

\newcommand{\bveps}{\mbox{\boldmath $\varepsilon$}}
\newcommand{\hSig}{\widehat\Sig}

\newcommand{\Sig}{\mathbf{\Sigma}}

\def\bI{{\bf I}}

\def\wcon{\stackrel{\cal D}{\longrightarrow}}

\def\bw{{\bf w}}
\def\be{{\bf e}}

\def\newpage{\vfill\eject}

\def\Var{\mbox{Var}}

\def\today{\ifcase\month\or
  January\or February\or March\or April\or May\or June\or
  July\or August\or September\or October\or November\or December\fi
  \space\number\day, \number\year}

\newdimen\biblioindent    \biblioindent=30pt

  at 7truept

\def\sgn{\mbox{sgn}}

\def\bI{{\bf I}}

\def\diag{\mbox{diag}}

\def\argmin{\mbox{argmin}}

\def\rss{(\bY-\bX\bbeta)^T(\bY-\bX\bbeta)}
\newcommand{\beq}{\begin{equation}}
\newcommand{\eeq}{\end{equation}}
\newcommand{\beqn}{\begin{eqnarray}}
\newcommand{\eeqn}{\end{eqnarray}}
\newcommand{\beqnn}{\begin{eqnarray*}}
\newcommand{\eeqnn}{\end{eqnarray*}}

\newcommand{\pln}{p_{\lambda_n}}
\newcommand{\plj}{p_{\lambda_j}}

\newcommand{\etal}{{\it et al }}

\newcommand{\gcv}{\mbox{GCV}}

\newcommand{\bZ}{\mbox{\bf Z}}

\newtheorem{theorem}{Theorem}

\theoremstyle{definition}

\def\rss{\mbox{RSS}}
\def\adj{\mbox{\footnotesize adj}}

\title[Statistical Challenges with High Dimensionality]
{Statistical Challenges with High Dimensionality: Feature Selection
in Knowledge Discovery}

\author[Jianqing Fan and Runze Li]{Jianqing Fan and Runze Li\thanks{Fan's research
was supported partially by NSF grant DMS-0354223, DMS-0532370 and
NIH R01-GM072611. Li's research was supported by NSF grant
DMS-0348869 and National Institute on Drug Abuse grant P50 DA10075.
The authors would like to thank Professors Peter Hall and Michael
Korosok for their constructive comments and John Dziak for his
assistance.}}

\begin{document}

\begin{abstract}
Technological innovations have revolutionized the process of
scientific research and knowledge discovery.  The availability of
massive data and challenges from frontiers of research and
development have reshaped statistical thinking, data analysis and
theoretical studies.  The challenges of high-dimensionality arise in
diverse fields of sciences and the humanities, ranging from
computational biology and health studies to financial engineering
and risk management.  In all of these fields, variable selection and
feature extraction are crucial for knowledge discovery. We first
give a comprehensive overview of statistical challenges with high
dimensionality in these diverse disciplines. We then approach the
problem of variable selection and feature extraction using a unified
framework: penalized likelihood methods. Issues relevant to the
choice of penalty functions are addressed. We demonstrate that for a
host of statistical problems, as long as the dimensionality is not
excessively large, we can estimate the model parameters as well as
if the best model is known in advance. The persistence property in
risk minimization is also addressed. The applicability of such a
theory and method to diverse statistical problems is demonstrated.
Other related problems with high-dimensionality are also discussed.
\end{abstract}



\maketitle

\section{Introduction}

Technological innovations have had deep impact on society and on
scientific research.  They allow us to collect massive amount of
data with relatively low cost.  Observations with curves, images or
movies, along with many other variables, are frequently seen in
contemporary scientific research and technological development. For
example, in biomedical studies, huge numbers of magnetic resonance
images (MRI) and functional MRI data are collected for each subject
with hundreds of subjects involved. Satellite imagery has been used
in natural resource discovery and agriculture, collecting thousands
of high resolution images. Examples of these kinds are plentiful in
computational biology, climatology, geology, neurology, health
science, economics, and finance among others. Frontiers of science,
engineering and the humanities differ in the problems of their
concerns, but nevertheless share one common theme: massive and
high-throughput data have been collected and new knowledge needs to
be discovered using these data.   These massive collections of data
along with many new scientific problems create golden opportunities
and significant challenges for the development of mathematical
sciences.

The availability of massive data along with new scientific problems
have reshaped statistical thinking and data analysis. Dimensionality
reduction and feature extraction play pivotal roles in all
high-dimensional mathematical problems.  The intensive computation
inherent in these problems has altered the course of methodological
development.  At the same time,  high-dimensionality has
significantly challenged traditional statistical theory.  Many new
insights need to be unveiled and many new phenomena need to be
discovered.  There is little doubt that the high dimensional data
analysis will be the most important research topic in statistics in
the 21st century \cite{Do00}.

Variable selection and feature extraction are fundamental to
knowledge discovery from massive data. Many variable selection
criteria have been proposed in the literature. Parsimonious models
are always desirable as they provide simple and interpretable
relations among scientific variables in addition to reducing
forecasting errors.  Traditional variable selection such as $C_p$,
AIC and BIC 
involves a combinatorial optimization problem, which is NP-hard,
with computational time increasing exponentially with the
dimensionality. The expensive computational cost makes traditional
procedures infeasible for high-dimensional data analysis. Clearly,
innovative variable selection procedures are needed to cope with
high-dimensionality.

Computational challenges from high-dimensional statistical endeavors
forge cross-fertilizations among applied and computational
mathematics, machine learning, and statistics.  For example, Donoho
and Elad \cite{DEl03} show that the NP-hard best subset regression
can be solved by a penalized $L_1$ least-squares problem, which can
be handled by a linear programming, when the solution is
sufficiently sparse. Wavelets are widely used in statistics function
estimation and signal processing \cite{AF01, CW92, Da92, DJ94,
DJKP95, Ma89a, Ma89b, Nik99}. Algebraic statistics, the term coined
by Pistone, Riccomagno, Wynn \cite{PRW00}, uses polynomial algebra
and combinatorial algorithms to solve computational problems in
experimental design and discrete probability \cite{PRW00},
conditional inferences based on Markovian chains \cite{DSt98},
parametric inference for biological sequence analysis \cite{PSt04},
and phylogenetic tree reconstruction \cite{SSu05}.

In high-dimensional data mining, it is helpful to distinguish two
types of statistical endeavors. In many machine learning problems
such as tumor classifications based on microarray or proteomics data
and asset allocations in finance, the interests often center around
the classification errors, or returns and risks of selected
portfolios 
rather than the accuracy of estimated parameters. On the other hand,
in many other statistical problems, concise relationship among
dependent and independent variables are needed. For example, in
health studies, we need not only to identify risk factors, but also
to assess accurately their risk contributions. These are needed for
prognosis and understanding the relative importance of risk factors.
Consistency results are inadequate for assessing the uncertainty in
parameter estimation. The distributions of selected and estimated
parameters are needed. Yet, despite extensive studies in classical
model selection techniques, no satisfactory solutions have yet been
produced.

In this article, we address the issues of variable selection and
feature extraction using a unified framework:  penalized likelihood
methods. This framework is applicable to both machine learning and
statistical inference problems.  In addition, it is applied to both
exact and approximate statistical modeling. We outline, in Section
2, some high-dimensional problems from computational biology,
biomedical studies, financial engineering, and machine learning, and
then provide a unified framework to address the issues of feature
selection in Sections 3 and 4.  In Sections 5 and 6, the framework
is then applied to provide solutions to some  problems outlined in
Section 2.

\section{Challenges from sciences and humanities}

 We now outline a few problems from various frontiers of research
to illustrate the challenges of high-dimensionality.  Some solutions
to these problems will be provided in Section 6.

\subsection{Computational biology}
Bioinformatic tools have been widely applied to genomics,
proteomics, gene networks, structure prediction, disease diagnosis
and drug design. The breakthroughs in biomedical imaging technology
allow scientists to monitor large amounts of diverse information on
genetic variation, gene and protein functions, interactions in
regulatory processes and biochemical pathways.   Such technology has also been
widely used for studying neuron activities and networks. Genomic
sequence analysis permits us to understand the homologies among
different species and infer their biological structures and
functionalities. Analysis of the network structure of protein can
predict the protein biological function. These quantitative
biological problems raise many new statistical and computational
problems. Let us focus specifically on the analysis of microarray
data to illustrate some challenges with dimensionality.

DNA microarrays have been widely used in simultaneously monitoring
mRNA expressions of thousands of genes in many areas of biomedical
research.  There are two popularly-used techniques: c-DNA
microarrays \cite{BB99} 
and Affymetrix GeneChip arrays \cite{LFG99}.
The former measures the abundance of mRNA expressions by mixing mRNAs
of treatment and control cells or tissues, hybridizing with cDNA on
the chip. The latter uses combined intensity information from
11-20 probes interrogating a part of the DNA sequence of a gene,
measuring separately mRNA expressions of treatment and control cells
or tissues.  Let us focus further on the cDNA microarray data.

The first statistical challenge is to remove systematic biases due
to experiment variations such as intensity effect in the scanning
process, block effect, dye effect, batch effect, amount of mRNA, DNA
concentration on arrays, among others.  This is collectively
referred to as normalization in the literature.  Normalization is
critical for multiple array comparisons.  Statistical models are
needed for estimation of these systematic biases in presence of
high-dimensional nuisance parameters from treatment effects on
genes. See, for example, lowess normalization in \cite{DYCS02, TOR01}, 
semiparametric model-based normalization
by \cite{FPH05,FTV04,HWZ05},
and robust normalization in \cite{MKH06}. 
The number of significantly expressed genes is relatively small.  Hence, model selection
techniques can be used to exploit the sparsity. In Section 6.1, we
briefly introduce semiparametric modeling techniques to issues of
normalization of cDNA microarray.

Once systematic biases have been removed, the statistical challenge
becomes selecting statistically significant genes based on a
relatively small sample size of arrays (e.g. $n=4, 6, 8$).  Various
testing procedures have been proposed in the literature.  See, for
example, \cite{FCC05, FTV04,HWZ05, TOR01,TTC01}.
In carrying out simultaneous testing
of orders of hundreds or thousands of genes, classical methods of controlling
the probability of making one falsely discovered gene are no longer
relevant.  Therefore various innovative methods have been proposed to control
the false discovery rates. See, for example, \cite{BY01,
DJ04,DSB03,Ef04,GW04,LRS05,STS04}.
The fundamental assumption in these developments is that the null
distribution of test statistics can be determined accurately.  This
assumption is usually not granted in practice and new probabilistic
challenge is to answer the questions how many simultaneous
hypotheses can be tested before the accuracy of approximations of
null distributions becomes poor. Large deviation theory
\cite{Hal87,Hal06, JSW03} is expected to play a critical role in
this endeavor. 
Some progress has been made using maximal inequalities \cite{KMa07}.

Tumor classification and clustering based on microarray and
proteomics data are another important class of challenging problems
in computational biology. Here, hundreds or thousands of gene
expressions are potential predictors, and the challenge is to select
important genes for effective disease classification and clustering.
See, for example, \cite{SND03,THN03,ZYS03}
for an overview and references therein.

Similar problems  include time-course microarray
experiments used to determine the expression pathways over time
\cite{SND03,TS06}
and genetic networks used for understanding interactions in
regulatory processes and biochemical pathways \cite{LG05}.
Challenges of selecting significant genes over time and classifying
patterns of gene expressions remain. In addition, understanding
genetic network problems requires estimating a huge covariance
matrix with some sparsity structure.
We introduce a modified Cholesky decomposition technique for
estimating large scale covariance matrices in Section 6.1.

\medskip

\subsection{Health studies}

Many health studies are longitudinal: each subject is followed over
a period of time and many covariates and responses of each subject
are collected at different time points. Framingham Heart Study
(FHS), initiated in 1948, is one of the most famous classic
longitudinal studies. Documentation of its first 50 years can be
found at the website of National Heart, Lung and Blood Institute
 (http://www.nhlbi.nih.gov/about/framingham/).
One can learn more details about this study
from the website of American Heart Association. In brief, the FHS
follows a representative sample of 5,209 adult residents and their
offspring aged 28-62 years in Framingham, Massachusetts. These
subjects have been tracked using (a) standardized biennial
cardiovascular examination, (b) daily surveillance of hospital
admissions, (c) death information and (d) information from
physicians and other sources outside the clinic.

In 1971 the study enrolled a second-generation group to participate
in similar examinations. It consisted of 5,124 of the original
participants' adult children and their spouses. This second study is
called the Framingham Offspring Study.

The main goal of this study is to identify major risk factors
associated with heart disease, stroke and other diseases, and to
learn the circumstances under which cardiovascular diseases arise,
evolve and end fatally in the general population. The findings in
this studies created a revolution in preventive medicine, and
forever changed the way the medical community and general public
view on the genesis of disease. In this study, there are more than
25,000 samples, each consisting of more than 100 variables. Because
of the nature of this longitudinal study, some participant cannot be
followed up due to their migrations. Thus, the collected data
contain many missing values. During the study, cardiovascular
diseases may develop for some participants, while other participants
may never experience with cardivoscular diseases. This implies that
some data are censored because the event of particular interest
never occurred. Furthermore, data between individuals may not be
independent because data for individuals in a family are clustered
and likely positively correlated. Missing, censoring and clustering
are common features in health studies. These three issues make data
structure complicated and identification of important risk factors
more challenging. In Section 6.2, we present a penalized partial
likelihood approach to selecting significant risk factors for
censored and clustering data. The penalized likelihood approach has
been used to analyze a data subset of Frammingham study in
\cite{CFLZ05}.

High-dimensionality is frequently seen in many other biomedical
studies. For example, ecological momentary assessment data have been
collected for smoking cessation studies.  In such a study, each
of a few hundreds participants
is provided a hand-held computer, which is designed to randomly
prompt the participants five to eight times per day over a period of
about 50 days and to provide 50 questions at each prompt. Therefore,
the data consist of a few hundreds of subjects and each of them may
have more than ten thousand observed values \cite{LRS06}. 
Such data are termed intensive longitudinal data. Classical
longitudinal methods are inadequate for such data. Walls and Schafer
\cite{WS06} presents more examples of intensive longitudinal data
and some useful models to analyze this kind of data.

\medskip

\subsection{Financial engineering and risk management}
Technological revolution and trade globalization have introduced a
new era of financial markets.  Over the last three decades, an
enormous number of new financial products have been created to meet
customers' demands.  For example, to reduce the impact of the
fluctuations of currency exchange rates on corporate finances, a
multinational corporation may decide to buy  options on the future
of exchange rates; to reduce the risk of price fluctuations of a
commodity (e.g. lumbers, corns, soybeans), a farmer may enter into a
future contract of the commodity; to reduce the risk of weather
exposures, amusement parks and energy companies may decide to
purchase financial derivatives based on the weather.  Since the
first options exchange opened in Chicago in 1973, the derivative
markets have experienced extraordinary growth. Professionals in
finance now routinely use sophisticated statistical techniques and
modern computing power in portfolio management, securities
regulation, proprietary trading, financial consulting, and risk
management. For an overview, see \cite{Fan05} 
and references therein.

Complex financial markets \cite{Hul03} make portfolio allocation,
asset pricing
and risk management very challenging. 
For example, the price of a stock depends not only on its past values,
but also its bond and derivative prices.  In addition, it depends on
prices of  related companies and their derivatives, and on overall
market conditions. Hence, the number of variables that influence
asset prices can be huge and the statistical challenge is to select
important factors that capture the market risks. Thanks to
technological innovations, high-frequency financial data are now
available for an array of different financial instruments over a
long time period. The amount of financial data available to financial
engineers is indeed astronomical.

Let us focus on a specific problem to illustrate the challenge of
dimensionality.  To optimize the performance of a portfolio
\cite{CLM97, Coc01}
or to manage the risk of a portfolio \cite{Mof03},
we need to estimate the covariance matrix of the
returns of assets in the portfolio. Suppose that we have 200 stocks
to be selected for asset allocation. There are 20,200 parameters in
the covariance matrix.  This is a high-dimensional statistical
problem and estimating it accurately poses challenges.

Covariance matrices pervade every facet of financial econometrics,
from asset allocation, asset pricing, and risk management, to derivative
pricing and proprietary trading.  As mentioned earlier, they are
also critical for studying genetic networks \cite{LG05}, 
as well as other statistical applications such as climatology
\cite{Joh01}.
In Section 6.1, a modified Cholesky decomposition is used to
estimate huge covariance matrices using penalized least squares
approach proposed in Section 2.
We will introduce a factor model for covariance estimation in Section
6.3.

\medskip

\subsection{Machine learning and data mining}

Machine learning and data mining extend traditional statistical
techniques to handle problems with much higher dimensionality.  The
size of data can also be astronomical:  from grocery sales and
financial market trading to biomedical images and natural resource
surveys. For an introduction, see the books \cite{HMS01, HTF01}.
Variable selections and feature extraction are vital for such
high-dimensional statistical explorations. Because of the size and
complexity of the problems, the associated mathematical theory also
differs from the traditional approach.  The dimensionality of
variables is comparable with the sample size and can even be much
higher than the sample size.  Selecting reliable predictors to
minimize risks of prediction is fundamental to machine learning and
data mining.  On the other hand, as the interest mainly lies in risk
minimization, unlike traditional statistics, the model parameters
are only of secondary interest. As a result, crude consistency
results suffice for understanding the performance of learning
theory.   This eases considerably the mathematical challenges of
high-dimensionality. For example, in the supervised (classification)
or unsupervised (clustering) learning, we do not need to know the
distributions of estimated coefficients in the underlying model.  We
only need to know the variables and their estimated parameters in
the model.
This differs from high-dimensional statistical problems in health
sciences and biomedical studies, where statistical inferences are
needed in presence of high-dimensionality. In Sections 4.2 and 6.4,
we will address further the challenges in machine learning.

\section{Penalized least squares}

With the above background, we now consider the variable selection in
the least-squares setting to gain further insights.  The idea will
be extended to the likelihood or pseudo-likelihood setting in the
next section. We demonstrate how to directly apply the penalized
least squares approach for function estimation or approximation
using wavelets or spline basis, based on noisy data in Section 5.
The penalized least squares method will be further extended to
penalized empirical risk minimization for machine learning in
Section 6.4.

Let $\{\bx_i, y_i\}, i=1,\cdots, n$, be a random sample from the
linear regression model
\beq y=\bx^T\bbeta+\eps, \label{e2.1}
\eeq
where $\eps$ is a random error with mean 0 and finite variance
$\sigma^2$, and $\bbeta=(\beta_1,\cdots,\beta_d)^T$ is the vector of
regression coefficients.
Here, we assume that all important
predictors, and their interactions or functions are already in the
model so that the full model (\ref{e2.1}) is correct. 

Many variable selection criteria or procedures are closely related
to minimize the following penalized least squares (PLS)
 \beq
  \frac 1{2n}\sum_{i=1}^n(y_i-\bx_i^T\bbeta)^2+\sum_{j=1}^dp_{\lambda_j}
  (|\beta_j|),
   \label{e2.2}
\eeq where $d$ is the dimension of $\bx$, and $p_{\lambda_j} (\cdot)$
is a penalty function, controlling model complexity.  The dependence
of the penalty function on $j$ allows us to incorporate prior
information. For instance, we may wish to keep certain important
predictors in the model and choose not to penalize their
coefficients.



The form of $\plj(\cdot)$ determines the general behavior of the
estimator. With the entropy or $L_0$-penalty, namely, $
\plj(|\beta_j|)=\frac 12\lambda^2I(|\beta_j|\ne0), $ the PLS
(\ref{e2.2}) becomes
 \beq
  \frac 1{2n}\sum_{i=1}^n (y_i-\bx_i^T\bbeta)^2 +
   \frac 12\lambda^2 |M|, \label{e2.3}
 \eeq
where $|M|=\sum_{j}I(|\beta_j|\ne0)$, the size of the candidate
model.  Among models with $m$ variables, the selected model is the
one with the minimum residual sum of squares (RRS), denoted by
$\rss_m$. A classical statistical method is to choose $m$ by
maximizing the adjusted $R^2$, given by
$$
      R_{\adj, m} = 1- \frac{n-1}{n-m}\frac{\rss_m}{\rss_1},
$$
or equivalently by minimizing $\rss_m/(n-m)$, where $\rss_1$ is the
total sum of squares based on the null model (using the intercept
only). Using $\log(1+x) \approx x $ for small $x$, it follows that
 \beq
   \log\{\rss_m /(n-m)\} 
    \approx (\log\sigma^2-1)+\sigma^{-2}\{\frac 1n \rss_m +
    \frac 1n m\sigma^2\}.\label{e2.4}
 \eeq
Therefore, maximization of $R_{\adj, m}$ is asymptotically equivalent to
minimizing  the PLS (\ref{e2.3}) with $\lambda = \sigma/\sqrt{n}$.
Similarly, generalized cross-validation (GCV) given by
\[
\gcv(m) = {\rss_m}/\{n(1-m/n)^2\}
\]
is asymptotically equivalent to the PLS (\ref{e2.3}) with
$\lambda=\sqrt{2}\sigma/\sqrt{n}$ and so is the cross-validation
(CV) criterion.

Many popular variable selection criteria can be shown asymptotically
equivalent to the PLS (\ref{e2.3}) with appropriate values of
$\lambda$, though these criteria were motivated from different
principles. See \cite{Mi02} 
 and references therein.
For instance, RIC \cite{FG94} 
corresponds to $\lambda
= \sqrt{2\log(d)}(\sigma/\sqrt{n})$. Since the entropy penalty
function is discontinuous, minimizing the entropy-penalized
least-squares requires exhaustive search, which is not feasible for
high-dimensional problem.  In addition, the sampling distributions
of resulting estimates are hard to derive.

Many researchers have been working on minimizing the PLS
(\ref{e2.2}) with $L_p$-penalty for some $p>0$. It is well known
that the $L_2$-penalty results in a ridge regression estimator,
which regularizes and stabilizes the estimator but introduces
biases. However, it does not shrink any coefficients directly to
zero.

The $L_p$-penalty with $0<p<2$ yields bridge regression \cite{FF93},
intermediating the best-subset ($L_0$-penalty) and the ridge
regression ($L_2$-penalty). The non-negative garrote \cite{Br95}
shares the same spirit as that of bridge regression. With the
$L_1$-penalty specifically, the PLS estimator is called
LASSO in \cite{Ti96}. 
In a seminal paper, Donoho and Elad \cite{DEl03} show that penalized
$L_0$-solution can be found by using penalized $L_1$-method for
sparse problem.  When $p\leq 1$, the PLS automatically performs
variable selection by removing predictors with very small estimated
coefficients.

Antoniadis and Fan \cite{AF01} discussed how to choose a penalty
function for wavelets regression. Fan and Li \cite{FL01}
advocated penalty
functions with three properties:

\begin{description}
\item[{a.  Sparsity:}] The resulting estimator should automatically set small
 estimated coefficients to zero to accomplish variable selection.

\item[{b. Unbiasedness:}] The resulting estimator should have low bias,
    especially when the true coefficient $\beta_{j}$ is large.

\item[{c. Continuity:}] The resulting estimator should be continuous
 to reduce instability in model prediction.
\end{description}

\noindent To gain further insights, let us assume that the design matrix
$\bX=(\bx_1,\cdots,\bx_n)^T$ for model (\ref{e2.1}) is orthogonal
and satisfies that $\frac 1n \bX^T\bX= I_d$. Let 
$\bz=(\bX^T\bX)^{-1} \bX^T\by$ be the least squares estimate of
$\bbeta$. Then (\ref{e2.2}) becomes
$$
\frac{1}{2n} \| \by - \bX \bz \| + \frac12 \|\bz-\bbeta \|^2  +
\sum_{j=1}^d \plj(|\beta_j|).
$$
Thus, the PLS reduces to a componentwise minimization problem:
\[
\min_{\beta_j} \{ \frac 12 (z_j-\beta_j)^2+\plj(|\beta_j|)\},
\quad\mbox{for} \quad j=1,\cdots, d,
\]
where $z_j$ is the $j$-th component of $\bz$.  Suppress the
subscript $j$ and  let
 \beq
   Q(\beta)=\frac12 (z-\beta)^2 +\pl(|\beta|). \label{e2.5}
 \eeq
Then the first order derivative of $Q(\beta)$ is given by
\[
Q'(\beta)=\beta-z+\pl'(|\beta|)\mbox{sgn}(\beta) =
\mbox{sgn}(\beta)\{|\beta|+\pl'(|\beta|)\} - z.
\]
Antoniadis and Fan \cite{AF01} and Fan and Li \cite{FL01} derived
that the PLS estimator possesses the following property:

(a) {\bf sparsity} if   $\min_{\beta} \{|\beta|+\pl'(|\beta|)\}>0$;

(b) {\bf unbiasedness}  $\pl'(|\beta|)=0$ for large $|\beta|$;

(c) {\bf continuity} if and only if $\mbox{argmin}_{\beta}\{|\beta|+\pl'(|\beta|)\}=0.$

\begin{figure}[htbp]  
\centerline{\psfig{figure=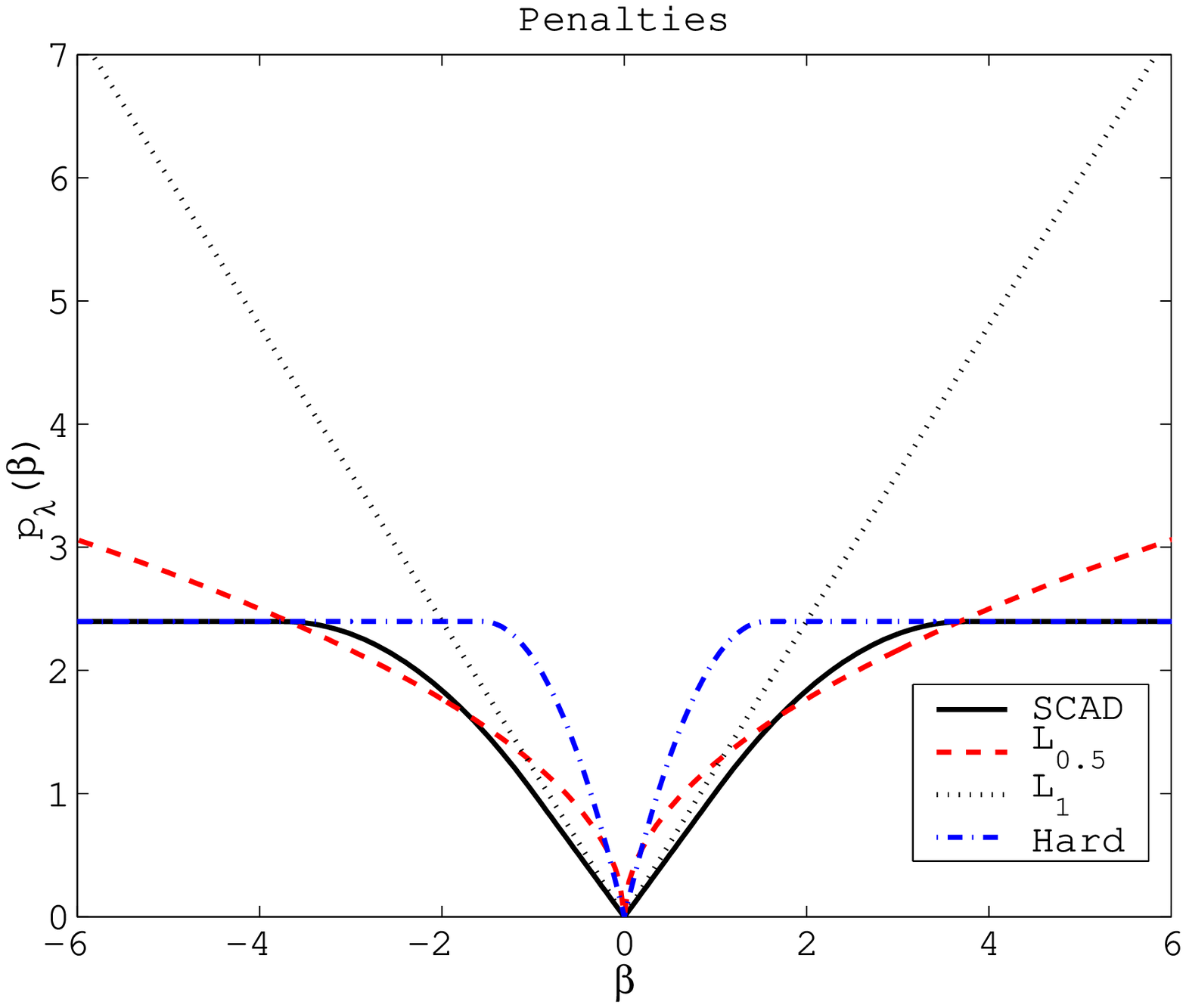,height=2.725in,width=2.3in}
\psfig{figure=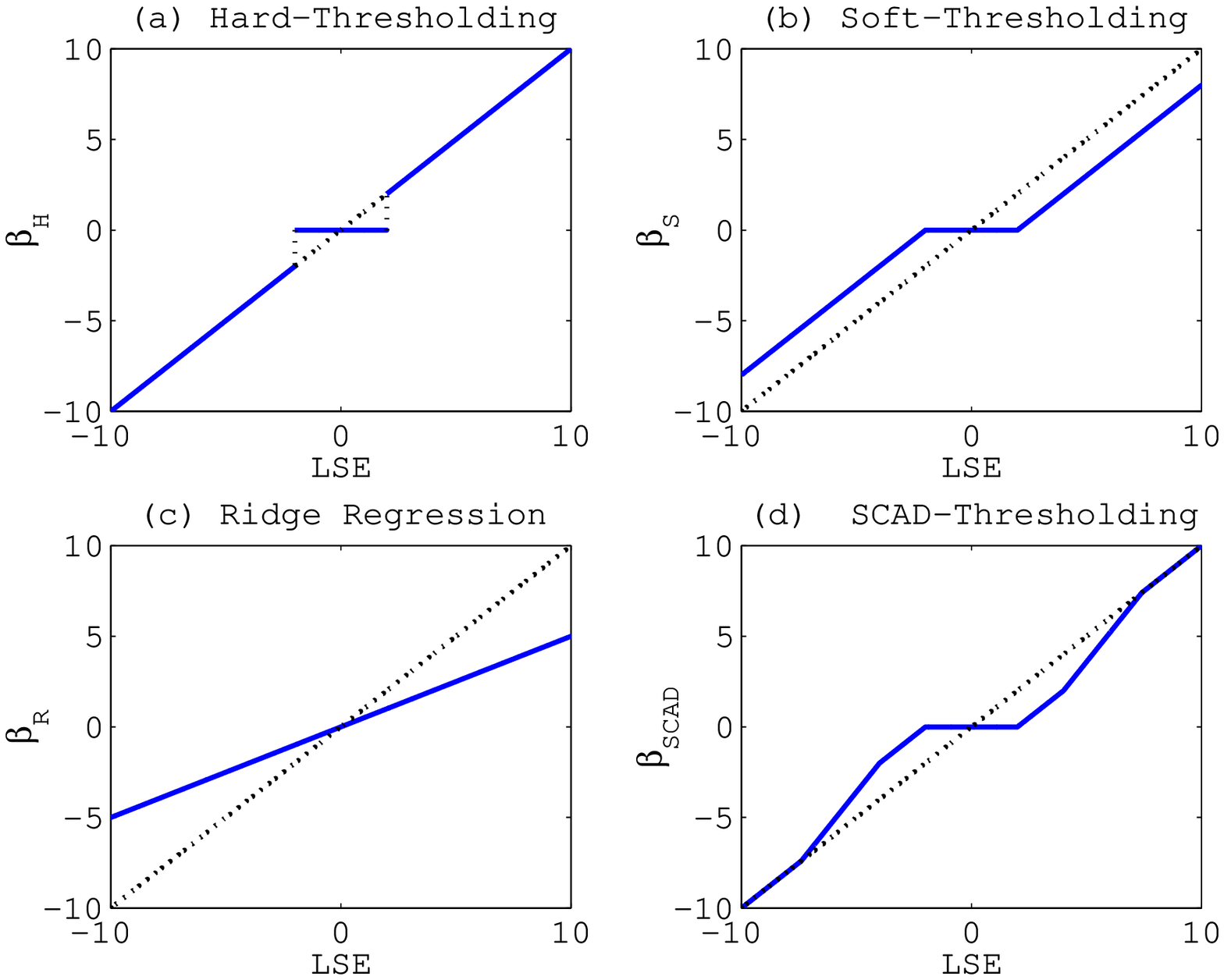,height=2.7in,width=2.7in} }
\caption{\label{FLfig1}   Penalty functions (left panel) and PLS
estimators (right panel).}
\end{figure}

\noindent The $L_p$-penalty with $0\le p<1$ does not satisfy the
continuity condition, the $L_1$ penalty does not satisfy the
unbiasedness condition, and $L_p$ with $p>1$ does not satisfy the
sparsity condition. Therefore, none of the $L_p$-penalties satisfies
the above three conditions simultaneously, and $L_1$-penalty is the
such penalty that is both convex and produces sparse solutions. Of
course, the class of penalty functions satisfying the aforementioned
three conditions are infinitely many. Fan and Li \cite{FL01}
suggested the use of the smoothly clipped absolute deviation (SCAD)
penalty defined as
\[
p_{\lambda}(|\beta|) = \left\{
\begin{array}{ll}
\lambda|\beta|, & \mbox{if} \ 0\le |\beta|<\lambda;\\
-(|\beta|^2 - 2 a \lambda |\beta| + \lambda^2)/\{2(a-1)\},\quad &
\mbox{if} \ \lambda
\le |\beta| < a\lambda;\\
(a+1)\lambda^2/2, & \mbox{if} \ |\beta|\ge a \lambda.
\end{array}
\right.
\]
They further suggested using $a=3.7$.  This function has similar
feature to the penalty function $\lambda |\beta|/(1+|\beta|)$
advocated in \cite{Nik99}.
Figure \ref{FLfig1} depicts the SCAD, $L_{0.5}$-penalty,
$L_1$-penalty, and hard thresholding penalty (to be introduced)
functions. These four penalty functions are singular at the origin,
a necessary condition for sparsity in variable selection.
Furthermore, the SCAD, hard-thresholding and $L_{0.5}$ penalties are
nonconvex over $(0,+\infty)$ in order to reduce the estimation bias.

Minimizing the PLS (\ref{e2.5}) with the entropy penalty or
hard-thresholding penalty $p_\lambda(\beta) = \lambda^2 - (\lambda -
|\beta|)_+^2$ (which is smoother) yields the hard-thresholding rule
\cite{DJ94}
$\hat{\beta}_H = z I(|z|>\lambda)$. With the $L_1$-penalty, the PLS
estimator is $ \hat{\beta}_S =  \sgn(z) (|z|-\lambda)_+$, the
soft-thresholding rule \cite{Bi83, DJ94}.
The $L_2$-penalty results in the ridge regression
$\hat{\beta}_R=(1+\lambda)^{-1}z$ and the SCAD penalty gives the
solution
\[
\hat{\beta}_{\mbox{\scriptsize SCAD}}= \left\{ \begin{array}{ll}
\sgn(z)(|z|-\lambda)_+,  &\ \mbox{when}\  |z|\le 2\lambda;\\
\{(a-1)z-\sgn(z)a\lambda\}/(a-2), &\ \mbox{when}\  2\lambda <|z|\le
a\lambda; \\ z, & \mbox{\ when \ } |z| >a\lambda.
\end{array} \right.
\]
These functions are also shown in Figure~\ref{FLfig1}.
The SCAD is an improvement over the $L_0$-penalty in two aspects:
saving computational cost and resulting in a continuous solution to
avoid unnecessary modeling variation. Furthermore, the SCAD improves
bridge regression by reducing modeling variation in model
prediction. Although similar in spirit to the $L_1$-penalty, the
SCAD also improves the  $L_1$-penalty by avoiding excessive
estimation bias since the solution of the $L_1$-penalty  could
shrink all regression coefficients by a constant, e.g., the soft
thresholding rule.



\section{Penalized likelihood}

PLS can easily be extended to handle a variety of response
variables, including binary response,  counts, and continuous
response. A popular family of this kind is called generalized linear
models. Our approach can also be applied to the case where the
likelihood is a quasi-likelihood or other discrepancy functions.
This will be demonstrated in Section 6.2 for analysis of survival
data, and in Section 6.4 for machine learning.

Suppose that conditioning on $\bx_i$, $y_i$ has a density
$f\{g(\bx_i^T\bbeta), y_i\}$, where $g$ is a known inverse link
function. Define a penalized likelihood as
\beq
    Q(\bbeta)=\frac 1n
    \sum_{i=1}^n \log f\{g(\bx_i^T\bbeta), y_i\} -  \sum_{j=1}^
    d\plj(|\beta_j|). \label{c1}
\eeq
Maximizing the penalized likelihood results in a penalized
likelihood estimator.
For certain penalties, such as the SCAD, the selected model based on
the nonconcave penalized likelihood satisfies $\beta_j=0$ for
certain $\beta_j$'s. Therefore, parameter estimation is performed at
the same time as the model selection.

\smallskip

\noindent{\bf Example 1}. ({\em Logistics Regression}) Suppose that
given $\bx_i$, $y_i$ follows a Bernoulli distribution with success
probability $P\{y_i=1|\bx_i\}=p(\bx_i)$. Take
$g(u)=\exp(u)/(1+\exp(u))$, i.e. $p(\bx) = \exp(\bx^T\bbeta)/\{1+
\exp(\bx^T\bbeta)\}$. Then (\ref{c1}) becomes
$$
\frac 1n \sum_{i=1}^n [y_i(\bx_i^T\bbeta) - \log\{1+\exp(\bx_i^T\bbeta)\}] -  \sum_{j=1}^
d\plj(|\beta_j|).
$$
Thus, variable selection for logistics regression can be achieved by
maximizing the above penalized likelihood.

\smallskip

\noindent{\bf Example 2}. ({\em Poisson Log-linear Regression}) Suppose that
given $\bx_i$, $y_i$ follows a Poisson distribution with mean
$\lambda(\bx_i)$. Take $g(\cdot)$ to be
the log-link, i.e. $\lambda(\bx) = \exp(\bx^T\bbeta)$. Then (\ref{c1}) can be written
as
$$
\frac 1n \sum_{i=1}^n \{y_i(\bx_i^T\bbeta) - \exp(\bx_i^T\bbeta)\} -  \sum_{j=1}^
d\plj(|\beta_j|)
$$
after dropping a constant. Thus, maximizing the above  penalized
likelihood with certain penalty functions yields a sparse solution
for $\bbeta$.

\subsection{Oracle properties}

Maximizing a penalized likelihood selects variables and estimates
parameters simultaneously.  This allows us to establish the sampling
properties of the resulting estimators. Under certain regularity
conditions, Fan and Li \cite{FL01} demonstrated how the rates of
convergence for the penalized likelihood estimators depend on the
regularization parameter $\lambda_n$ and established the oracle
properties of the penalized likelihood estimators.

In the context of variable selection for high-dimensional modeling,
it is natural to allow the number of introduced variables to grow
with the sample sizes. Fan and Peng \cite{FP04} have studied the
asymptotic properties of the penalized likelihood estimator for
situations in which the number of parameters, denoted by $d_n$,
tends to $\infty$ as the sample size $n$ increases. Denote
$\bbeta_{n0}$ to be the true value of $\bbeta$. To emphasize the
dependence of $\lambda_j$ on $n$, we use notation $\lambda_{n,j}$
for $\lambda_j$ in this subsection. Define \beq
 a_n = \max \{ p_{\lambda_{n,j}}'(|\beta_{n0j}|): \beta_{n0j} \not = 0\}
\quad \mbox{and}\quad b_n= \max \{
|p_{\lambda_{n,j}}''(|\beta_{n0j}|)|: \beta_{n0j} \not = 0\}.
\label{c2} \eeq

Fan and Peng \cite{FP04} showed that if both $a_n$ and $b_n$ tend to
0 as $n\to\infty$, then under certain regularity conditions, there
exists a local maximizer $\hat{\bbeta}$ of $Q(\bbeta)$ such that
\beq \| \hat{\bbeta} - \bbeta_{n0} \| = O_P\{\sqrt{d_n}(n^{-1/2}
+a_n)\}. \label{c2b} \eeq
It is clear from (\ref{c2b}) that by choosing a proper
$\lambda_{n,j}$ such that $a_n = O(n^{-1/2})$, there exists a
root-$(n/d_n)$ consistent penalized likelihood estimator. For
example, for the SCAD, the penalized likelihood estimator is
root-$(n/d_n)$ consistent if all $\lambda_{n,j}$'s tend to 0.

Without loss of generality assume that, unknown to us, the first
$s_n$ components of $\bbeta_{n0}$, denoted by $\bbeta_{n01}$, are
nonzero and do not vanish and the remaining $d_n-s_n$ coefficients,
denoted by $\bbeta_{n02}$, are 0. Denote by
\[
\Sig=\mbox{diag}\left\{p''_{\lambda_{n,1}}(|\beta_{n01}|),\cdots,
p''_{\lambda_{n,s_n}}(|\beta_{n0s_n}|)\right\}
\]
and
\[
\bb
=\left(p'_{\lambda_{n,1}}(|\beta_{n01}|)\sgn(\beta_{n01}),
\cdots,p'_{\lambda_{n,s_n}}(|\beta_{n0s_n}|)\sgn(\beta_{n0s_n})\right)^T.
\]

\begin{theorem}  
Assume that as $n\to\infty$,
$
  \min_{1\le j\le s_n} |\beta_{n0j}|/\lambda_{n,j} \to \infty
 $
and that the penalty function $\plj(|\beta_j|)$ satisfies
\beq
  \liminf_{n\to\infty}\liminf_{\beta_j\to0+} p'_{\lambda_{n,j}}(\beta_j) /
   \lambda_{n,j} > 0. \label{c3}
\eeq If $\lambda_{n,j}\to 0$, $\sqrt{n/d_n}\lambda_{n,j}\to \infty$
and $d_n^5/n\to 0$ as $n\to\infty$, then with probability tending to
1, the root $n/d_n$  consistent local maximizers $
\hat{\bbeta}=(\hat{\bbeta}_{n1}^T, \hat{\bbeta}_{n2}^T)^T$ must
satisfy:
\begin{description}
\item{(i)} ({\bf Sparsity}) $\hat{\bbeta}_{n2}={\bf 0}$;
\item{(ii)} ({\bf Asymptotic normality}) for any $q\times s_n$ matrix
$\bA_n$ such that $\bA_n \bA_n^T\to \bG$, a $q\times q$ positive
definite symmetric matrix,
\[
\sqrt{n}\bA_n \bI_1^{-1/2}\{\bI_1+\Sig\} \left\{\hat{\bbeta}_{n1}-
\bbeta_{n10} + (\bI_1+\Sig)^{-1}\bb\right\} \wcon N({\bf 0}, \bG)
\]
where $\bI_1=\bI_1(\bbeta_{n10},{\bf 0})$, the Fisher information
knowing $\bbeta_{n20}={\bf 0}$.
\end{description}
\end{theorem}

The theorem implies that any finite set of elements of $\hat{\bbeta}_{n1}$
are jointly asymptotically normal. For the SCAD, if all
$\lambda_{j,n}\to 0$, $a_n = 0$. Hence, when
$\sqrt{n/d_n}\lambda_{n,j}\to \infty$, its corresponding penalized
likelihood estimators possess the oracle property, i.e.,  perform as
well as the maximum likelihood estimates for estimating
$\bbeta_{n1}$  knowing $\bbeta_{n2}={\bf 0}$. That is, with probability approaching to
1,
\[
 \hbbeta_{n2}=0,\quad\mbox{and} \quad
 \sqrt{n}\bA_n\bI_1^{1/2}(\hat{\bbeta}_{n1}- \bbeta_{n10})\to N({\bf 0},
 \bG).
\]

For the $L_1$-penalty, $a_n = \max_j\lambda_{j,n}$.  Hence, the
root-$n/d_n$ consistency requires that
$\lambda_{n,j}=O(\sqrt{d_n/n})$. On the other hand, the oracle
property in Theorem 2 requires that $\sqrt{n/d_n}\lambda_{n,j}\to
\infty$. These two conditions for LASSO cannot be satisfied
simultaneously. It has indeed been shown that the oracle property
does not hold for the  $L_1$-penalty even in the finite parameter
setting \cite{Zo05}. 

\subsection{Risk minimization and persistence}

In machine learning such as tumor classifications, the primary
interest centers on the misclassification errors or more generally
expected losses, not the accuracy of estimated parameters. This kind
of properties is called persistence in \cite{Gre06,GRi04}.

Consider predicting the response $Y$ using a class of model $g(\bx^T
\bbeta)$ with a loss function $\ell \{g(\bX^T \bbeta), Y)$. Then,
the risk is
$$
    L_n(\bbeta) = E \ell \{g(\bX^T \bbeta), Y\},
$$
where $n$ is used to stress the dependence of dimensionality $d$ on
$n$. The minimum risk is obtained at $\bbeta_n^* = \argmin_{\bbeta}
L_n(\bbeta)$.  In the penalized likelihood context, $\ell = - \log
f$. Suppose that there is an estimator $\hat{\bbeta}_n$ based on a
sample of size $n$.  This can be done by the penalized empirical
risk minimization similarly to (\ref{c1}):
\begin{equation}
  n^{-1} \sum_{i=1}^n \ell\{ g(\bx_i^T \bbeta), y_i \}
  +  \sum_{j=1}^d \plj(|\beta_j|),  \label{aaa}
\end{equation}
based on a set of training data $\{(\bx_i,y_i), i=1,\cdots,n\}$. The
persistence requires
\begin{equation}
     L_n(\hbbeta) - L_n (\bbeta^*) \stackrel{P}{\longrightarrow}      0,
\end{equation}
but not the consistency of $\hbbeta$ to $\bbeta_n^*$.  This is in
general a much weaker mathematical requirement. Greenshtein and
Ritov \cite{GRi04} show that if the non-sparsity rate $s_n =
O\{(n/\log n)^{1/2}\}$ and $d_n = n^{\alpha}$ for some $\alpha > 1$,
LASSO (penalized $L_1$ least-squares) is persistent under the
quadratic loss. Greenshtein \cite{Gre06} extends the results to the
case where $s_n = O\{n/\log n\}$ and more general loss functions.
Meinshausen \cite{Mei05} considers a case with finite non-sparsity
$s_n$ but with $\log d_n = n^\xi$, with $\xi \in (0, 1)$. It is
shown there that for the quadratic loss, LASSO is persistent, but
the rate to persistency is slower than a relaxed LASSO.  This again
shows the bias problems in LASSO.

\subsection{Issues in practical implementation}

In this section, we address practical implementation issues related
to the PLS and penalized likelihood.

\smallskip

\noindent{\bf Local quadratic approximation (LQA).} The $L_p$, $(0 <
p < 1$), and SCAD penalty functions are singular at the origin, and
they do not have continuous second order derivatives. Therefore,
maximizing the nonconcave penalized likelihood is challenging.
Fan and Li \cite{FL01} propose  locally approximating them by a
quadratic function as follows. Suppose that we are given an initial
value ${\mbox{\boldmath$\beta$}}^0$  that is close to the optimizer
of $Q(\bbeta)$. For example, take initial value to be the maximum
likelihood estimate (without penalty). Under some regularity
conditions, the initial value is  a consistent estimate for
$\bbeta$, and therefore it is close to the true value. Thus, we can
locally approximate the penalty function by a quadratic function as
\begin{equation}
p_{\lambda_n}(|\beta_j|) \approx p_{\lambda_n}(|\beta_{j}^0|) + \frac{1}{2}
    \{p_{\lambda_n}'(|\beta_{j}^0|)/ |\beta_{j}^0|\}
    (\beta_j^2 - \beta_{j}^{02}), \mbox{ for $\beta_j \approx
    \beta_{j}^0$}.
    \label{e1}
\end{equation}
To avoid numerical instability, we set $\hat{\beta}_j=0$ if
$\beta_{j}^0$ is very close to 0. This corresponds to deleting $x_j$
from the final model. With the aid of the LQA, the optimization of
penalized least-squares, penalized likelihood or penalized partial
likelihood (see Section 6.2) can be carried out by using the
Newton-Raphson algorithm. It is worth noting that the LQA should be
updated at each step during the course of iteration of the
algorithm. We refer to the modified Newton-Raphson algorithm as the
LQA algorithm.

The convergence property of the LQA algorithm was studied in
\cite{HL05}, whose authors first showed that the LQA plays the same role as
the E-step in the EM
algorithm \cite{DLR77}. 
Therefore the behavior of the LQA algorithm is similar to the EM
algorithm. Unlike the original EM algorithm, in which a full
iteration for maximization is carried out after every E-step, we
update the LQA at each step during the iteration course. This speeds
up the convergence of the algorithm. The convergence rate of the LQA
algorithm is quadratic which is the same as that of the modified EM
algorithm \cite{La95}.

When the algorithm converges,
the estimator satisfies the condition
\[
 {\partial \ell (\hat{{\mbox{\boldmath$\beta$}}})}/{\partial \beta_j} + n
 p_{\lambda_j}'(|\hat{\beta}_{j}|)
 \mbox{sgn}(\hat{\beta}_{j}) = 0,
\]
the penalized likelihood equation, for non-zero elements of
$\hat{{\mbox{\boldmath$\beta$}}}$. \medskip

\noindent{\bf Standard error formula.} Following conventional
techniques in the likelihood setting, we can estimate the standard
error of the resulting estimator by using the sandwich formula.
Specifically, the corresponding sandwich formula can be used as an
estimator for the covariance of the estimator
$\hat{{\mbox{\boldmath$\beta$}}}_1$, the non-vanishing component of
$\hat{{\mbox{\boldmath$\beta$}}}$.  That is,
\begin{equation}
\widehat{\mbox{cov}}(\hat{{\mbox{\boldmath$\beta$}}}_1)= \{\nabla^2
\ell(\hat{\mbox{\boldmath$\beta$}}_1) - n
\Sig_{\lambda}(\hat{\mbox{\boldmath$\beta$}}_1)\}^{-1}
\widehat{\mbox{cov}}\{\nabla \ell(\hat{\mbox{\boldmath$\beta$}}_1)\}
\{\nabla^2 \ell(\hat{\mbox{\boldmath$\beta$}}_1)- n
\Sig_{\lambda}(\hat{\mbox{\boldmath$\beta$}}_1)\}^{-1},
\end{equation}
where $\widehat{\mbox{cov}}\{\nabla
\ell(\hat{\mbox{\boldmath$\beta$}}_1)\}$ is the usual empirically estimated
covariance matrix and
$$
\Sig_{\lambda}(\hat{\mbox{\boldmath$\beta$}}_1)=\mbox{diag}\{
p'_{\lambda_1}(|\hat{\beta}_{1}|)/|\hat{\beta}_{1}|, \cdots,
p'_{\lambda_{s_n}}(|\hat{\beta}_{s_n}|)/|\hat{\beta}_{s_n}|\}
$$
and $s_n$ the dimension of $\hat{\mbox{\boldmath$\beta$}}_1$. Fan
and Peng \cite{FP04} 
demonstrated the consistency of the sandwich formula:

\begin{theorem}  
Under the conditions of Theorem 1, we have
$$
\bA_n \widehat{\mbox{cov}}(\hat{{\mbox{\boldmath$\beta$}}}_1)
\bA^T_n - \bA_n \Sig_n \bA^T_n \stackrel{P}{\longrightarrow} 0 \
\mbox{as} \ n \to \infty
$$
for any 
matrix $\bA_n$ such that $\bA_n\bA^T_n=\bG$, where $\Sig_n = (\bI_1
+ \Sig)^{-1} \bI_1^{-1} (\bI_1 + \Sig)^{-1}$.
\end{theorem}

\noindent{\bf Selection of regularization parameters.} To implement
the methods described in previous sections, it is desirable to have
an automatic method for selecting the thresholding parameter
$\lambda$ in $p_{\lambda}(\cdot)$ based on data. Here, we estimate
$\lambda$ via minimizing an approximate generalized cross-validation
(GCV) statistic in \cite{CW79}.
By some
straightforward calculation, the effective number of parameters for
$Q({\mbox{\boldmath$\beta$}})$ in the last step of the
Newton-Raphson algorithm iteration is
$$
e(\bla)\equiv e(\lambda_1,\cdots,\lambda_d)=\mbox{tr} [\{\nabla^2
\ell(\hat{\mbox{\boldmath$\beta$}})+
\Sig_{\lambda}(\hat{\mbox{\boldmath$\beta$}})\}^{-1}\nabla^2
\ell(\hat{\mbox{\boldmath$\beta$}}) ].
$$
Therefore the generalized cross-validation statistic is defined by
$$
\mbox{GCV}(\bla)=
{-\ell(\hat{\mbox{\boldmath$\beta$}})}/{[n\{1-e(\bla)/n\}^2]}
$$
and $\hat{\bla}=\mbox{argmin}_{\bla} \{\mbox{GCV}(\bla)\}$ is
selected.

To find an optimal $\bla$, we need to minimize the GCV over a
$d_n$-dimensional space. This is an unduly onerous task.
Intuitively, it is expected that the magnitude of $\lambda_j$ should
be proportional to the standard error of the maximum likelihood
estimate of $\beta_j$. Thus, we set $\bla = \lambda
\mbox{se}(\hbbeta_{\mbox{\footnotesize MLE}})$ in practice, where
$\mbox{se}(\hbbeta_{\mbox{\footnotesize MLE}})$ denotes the standard
error of the MLE. Therefore, we minimize the GCV score over the
one-dimensional space, which will save a great deal of computational
cost.  The behavior of such a method has been investigated recently.

\section{Applications to function estimation}



Let us begin with one-dimensional function estimation. Suppose that
we have noisy data at possibly irregular design points
$\{x_1,\cdots, x_n\}$:
\[
y_i=m(x_i) +\varepsilon_i,
\]
where $m$ is an unknown regression and $\eps_i$'s are iid random
error following $N(0,\sigma^2)$. Local modeling techniques
\cite{FG96} have been widely used to estimate $m(\cdot)$. Here we
focus on global function approximation methods.

Wavelet transforms are a device for representing functions in a way
that is local in both time and frequency domains
\cite{CS94,CW92,Ma89a,Ma89b}. During the last decade, they have
received a great deal of attention in applied mathematics, image
analysis, signal compression, and many other fields of engineering.
Daubechies \cite{Da92} and Meyer \cite{Me90} are good introductory
references to this subject.
Wavelet-based methods have many exciting statistical
properties \cite{DJ94}. 
Earlier papers on wavelets assume the regular design points, i.e,
$x_i=\frac i n$ (usually $n=2^k$ for some integer $k$) so that fast
computation algorithms can be implemented. See \cite{DJKP95}
and references therein.  For an overview of wavelets in statistics,
see \cite{Wan06}.

Antoniadis and Fan \cite{AF01} discussed how to apply wavelet
methods for function estimation with irregular design points using
penalized least squares. Without loss of generality, assume that
$m(x)$ is defined on $[0,1]$.  By moving nondyadic points to dyadic
points, we assume $x_i=n_i/2^J$ for some $n_i$ and some fine
resolution $J$ that is determined by  users. To make this
approximation errors negligible, we take $J$ large enough such that
$2^J\ge n$. Let $\bW$ be a given wavelet transform at all dyadic
points $\{i/2^J:, i=1,\cdots, 2^J-1\}$. Let $N=2^J$ and $\ba_i$ be
the $n_i$-th column of $\bW$, an $N\times N$ matrix, and $\bbeta =
\bW \mbox{\bf m}$ be the wavelet transform of the function $m$ at
dyadic points. Then, it is easy to see that $m(x_i)=\ba_i^T\bbeta$.
This yields an overparameterized linear model
 \beq
   y_i=\ba_i^T \bbeta + \eps_i, \label{g1}
 \eeq
which aims at reducing modeling biases. However, one cannot find a
reasonable estimate of $\bbeta$ by using the ordinary least squares
method since $N\ge n$. Directly applying penalized least squares, we
have
 \beq
   \frac {1}{2n} \sum_{i=1}^n ( y_i - \ba_i^T\bbeta)^2 + \sum_{j=1}^N
\plj(|\beta_j|). \label{g2} \eeq If the sampling points are equally
spaced and $n=2^J$, the corresponding design matrix of linear model
(\ref{g1}) becomes a square orthogonal matrix. From the discussion
in Section 3, minimizing the PLS (\ref{g2}) with the entropy penalty
or the hard-thresholding penalty results in a hard-thresholding
rule. With the $L_1$ penalty, the PLS estimator is the
soft-thresholding rule. Assume that $\pl(\cdot)$ is nonnegative,
nondecreasing, and differentiable over $(0,\infty)$ and that
function $-\beta-\pl'(\beta)$ is strictly unimodal on $(0,\infty)$,
$\pl'(\cdot)$ is nonincreasing and $\pl'(0+)>0$. Then, Antoniadis
and Fan \cite{AF01} showed that the resulting penalized
least-squares estimator that minimizes (\ref{g2}) is adaptively
minimax within a factor of logarithmic order as follows. Define the
Besov space ball $B_{p,q}^r(C)$ to be
$$
B_{p,q}^r(C) = \{m\in L_p: \sum_j
(2^{j(r+1/2-1/p)}\|\btheta_{j\cdot}\|_p)^q < C\},
$$
where $\btheta_{j\cdot}$ is the vector of wavelet coefficients of
function $m$ at the resolution level $j$. Here $r$ indicates the
degree of smoothness of the regression functions $m$.

\begin{theorem}  
Suppose that the regression function $m(\cdot)$ is in a Besov ball
with $r+1/2-1/p>0$. Then, the  maximum risk of the PLS estimator
$\hat{m}(\cdot)$ over $B_{p,q}^r(C)$ is of rate
$O(n^{-2r/(2r+1)}\log(n))$ when the universal thresholding
$\sqrt{2\log(n)/n}$ is used. It also achieves the rate of
convergence $O\{n^{-2r/(2r+1)}\log(n)\}$ when the minimax
thresholding $p_n/\sqrt{n}$ is used, where $p_n$ is given in
\cite{AF01}.
\end{theorem}

We next consider multivariate regression function estimation.
Suppose that $\{\bx_i,y_i\}$ is a random sample from the regression
model
\[
y=m(\bx) + \eps,
\]
where, without loss of generality, it is assumed that $\bx\in
[0,1]^d$. Radial basis and neural-network are also popular for
approximating multi-dimensional functions. In the literature of
spline smoothing, it is typically assumed that the mean function
$m(\bx)$ has a low-dimensional structure. For example,
\[
m(\bx) = \mu_0 + \sum_j m_j(x_j) + \sum_{k< l} m_{kl}(x_k,x_l).
\]
For given knots, a set of spline basis functions can be constructed.
The two most popular spline bases are the  truncated power spline basis
$1,x,x^2,x^3, (x-t_j)_+^3,\, (j=1,\,\cdots,\,J), $ where $t_j$'s are
knots, and the B-spline basis (see \cite{Bo78} 
for definition). The
$B$-spline basis is numerically more stable since the multiple
correlation among the basis functions is smaller, but the power
truncated  spline basis has the advantage that deleting a basis
function is the same as deleting a knot.

For a given set of 1-dimensional spline bases, we can further
construct a multivariate spline basis using tensor products. Let
$\{B_1,\,\cdots, B_{J}\}$ be a set of spline basis functions on
$[0,1]^d$.  Approximate the regression function $m(\bx)$ by a linear
combination of the basis functions, $\sum \beta_jB_j(\bx)$, say. To
avoid a large approximation bias,  we take a large $J$. This yields
an overparameterized linear model, and the fitted curve of the least
squares estimate is typically undersmooth. Smoothing
spline 
suggested penalizing the roughness of the
resulting estimate. This is equivalent to the penalized least
squares with a quadratic penalty. In a series of work by Stone and
his collaborators (see \cite{SHKT97}), 
they advocate using regression splines and modifying traditional
variable selection approaches to select useful spline subbases.
Ruppert  \etal \cite{RWC03} advocated penalized splines in
statistical modeling, in which power truncated splines are used with
the $L_2$ penalty. Another kind of penalized splines method proposed
by \cite{EM96}
shares the same spirit of \cite{RWC03}.

\section{Some solutions to the challenges}

In this section, we provide some solutions to problems raised in
Section 2.

\subsection{Computational biology}

As discussed in Section 2.1, the first statistical challenge in
computational biology is how to remove systematic biases due to
experiment variations. Thus, let us first discuss the issue of
normalization of cDNA-microarrays. Let $Y_g$ be the log-ratio of the
intensity of gene $g$ of the treatment sample over that of the
control sample.  Denote by $X_g$ the average of the log-intensities
of gene $g$ at the treatment and control samples. Set $r_g$ and
$c_g$ be the row and column of the block where the cDNA of gene $g$
resides. Fan \etal \cite{FTV04} use the following model to estimate
the intensity and block effect:
\begin{equation}
 Y_g = \alpha_g + \beta_{r_g} + \gamma_{c_g} + f(X_g) +
 \varepsilon_g, \quad g = 1, \cdots, N\label{f1}
\end{equation}
where $\alpha_g$ is the treatment effect on gene $g$, $\beta_{r_g}$
and $\gamma_{c_g}$ are block effects that are decomposed into the
column and row effect, $f(X_g)$ represents the intensity effect and
$N$ is the total number of genes. Based on $J$ arrays, an aim of
microarray data analysis is to find genes $g$ with $\alpha_g$
statistically significantly different from 0. However, before
carrying multiple array comparisons, the block and treatment effects
should first be estimated and removed. For this normalization
purpose, parameters $\alpha_g$ are nuisance and high-dimensional
(recall $N$ is in the order of tens of thousands). On the other
hand, the number of significantly expressed genes is relatively
small, yielding the sparsity structure of $\alpha_g$.

Model (\ref{f1}) is not identifiable.  Fan \etal \cite{FTV04} use
within-array replicates to infer about the block and treatment
effects. Suppose that we have $I$ replications for $G$ genes, which
could be a small fraction of $N$.  For example, in \cite{FTV04},
only 111 genes were repeated at random blocks ($G=111, I = 2$),
whereas in \cite{MKH06},
all genes were repeated three times, i.e. $I =3$ and  $N = 3 G$,
though both have about $N \approx 20,000$ genes printed on an array.
Using $I$ replicated data on $G$ genes, model (\ref{f1}) becomes
\begin{eqnarray}
Y_{gi} = \alpha_{g} + \beta_{r_{gi}} + \gamma_{c_{gi}} + f(X_{gi}) +
 \varepsilon_{gi},  \quad g = 1, \cdots, G; i=1, \cdots, I.\label{f2}
\end{eqnarray}
With estimated coefficients $\hat{\beta}$ and $\hat{\gamma}$ and the
function $\hat{f}$, model (\ref{f1}) implies that the normalized
data are $ Y_{g}^* = Y_g  - \hat{\beta}_{r_g} - \hat{\gamma}_{c_g} -
\hat{ f}(X_g)$ even for non-repeated genes.

Model (\ref{f2}) can be used to remove the intensity effect array by
array, though the number of nuisance parameters is very large, a
fraction of total sample size in (\ref{f2}).  To improve the
efficiency of estimation, Fan \etal \cite{FPH05} aggregate the
information from other microarrays (total $J$ arrays):
\begin{eqnarray}
  Y_{gij} = \alpha_{g} + \beta_{r_{gi}, j} + \gamma_{c_{gi}, j} + f_j(X_{gij}) +
   \varepsilon_{gi},  \quad j=1, \cdots, J,\label{f3}
\end{eqnarray}
where the subscript $j$ denotes the array effect.

The parameters in (\ref{f2}) can be estimated by the profile
least-squares method using the Gauss-Seidel type of algorithm.  See
\cite{FPH05} 
for details.  To state the results, let us write
model (\ref{f2}) as
\begin{eqnarray}
    Y_{gi} = \alpha_g + \bZ_{gi}^T \bbeta + f(X_{gi}) +
    \varepsilon_{gi}, \label{f4}
\end{eqnarray}
by appropriately introducing the dummy variable $\bZ$.  Fan \etal
\cite{FPH05} obtained the following results.

\begin{theorem}  
Under some regularity conditions, as $n = IG \to \infty$, the
profile least-squares estimator of model (\ref{f4}) has
$$
\sqrt{n} (\hbbeta - \bbeta) \wcon N \left (0, \frac{I}{I-1} \sigma^2
\Sig^{-1}  \right ),
$$
where $\Sig =  E\{ \Var(\bZ|X)\}$ and $\sigma^2 =
\Var(\varepsilon)$.  In addition, $\hat{f}(x) - f(x) =
O_P(n^{-2/5})$.
\end{theorem}

\begin{theorem}  
Under some regularity conditions, as $n = IG \to \infty$, when $X$
and $\bZ$ are independent, the profile least-squares estimator based
on (\ref{f3}) possesses
$$
\sqrt{n} (\hbbeta _j- \bbeta_j) \wcon N \left (0,
\frac{I(J-1)+1}{J(I-1)}\; \sigma^2 \Sig^{-1} \right ).
$$
\end{theorem}

The above theorems show that the block effect can be estimated at
rate $O_P(n^{-1/2})$ and intensity effect $f$ can be estimated at
rate $O_P(n^{-2/5})$.  This rate can be improved to $O_P(n^{-1/2} +
N^{-2/5})$ when data in (\ref{f1}) are all used.  The techniques
have also been adapted for the
normalization of Affymetrix arrays \cite{FCC05}. 
Once the
arrays have been normalized, the problem becomes selecting
significantly expressed genes using the normalized data
\begin{equation}
    Y_{gj}^* = \alpha_g + \varepsilon_{gj}, \qquad g =1, \cdots, N; j =1, \cdots,
    J,
\end{equation}
where $Y_{gj}^*$ is the normalized expression of gene $g$ in array
$j$. This is again a high-dimensional statistical inference problem.
The issues of computing P-values and false discovery are given in
Section  2.1.

Estimation of high-dimensional covariance matrices is critical in
studying genetic networks. PLS and penalized likelihood can be used
to estimate large scale covariance matrices effectively and
parsimoniously \cite{HLPL06, LD06}.
Let $\bw=(W_1,\cdots, W_d)^T$ be a $d$-dimensional random vector
with mean zero and covariance $\Sigma.$ Using the modified Cholesky
decomposition, we have $\bL\Sig \bL^T = \bD$, where $\bL$ is a lower
triangular matrix having ones on its diagonal and typical element
$-\phi_{tj}$ in the $(t,j)$th position for $1\le j< t\le d$, and
$\bD=\diag\{\sigma_1^2,\cdots,\sigma_d^2)^T$ is a diagonal matrix.
Denote  $\be=\bL\bw=(e_1,\cdots, e_d)^T$. Since $\bD$ is diagonal,
$e_1,\cdots, e_d$ are uncorrelated. Thus, for $2\le t\le d$
\beq
W_t=\sum_{j=1}^{t-1}\phi_{tj}W_{j} + e_t. \label{e6.2}
\eeq
That is,
the $W_t$ is an autoregressive (AR) series, which gives an
interpretation for elements of $L$ and $D$, and allows us to use PLS
for covariance selection. We  first estimate $\sigma_t^2$ using
the mean squared errors of model (\ref{e6.2}). Suppose that $\bw_i$, $i=1,\cdots,n$, is a
random sample from $\bw$.  For $t=2,\cdots,d$, covariance selection
can be achieved by minimizing the following PLS functions: \beq
    \frac 1{2n}\sum_{i=1}^n (W_{it}-
    \sum_{j=1}^{t-1}\phi_{tj}W_{ij})^2 + \sum_{j=1}^{t-1}
    p_{\lambda_{t,j}}(|\phi_{tj}|). \label{e6.4}
\eeq This reduces the non-sparse elements in the lower triangle
matrix $\bL$.  With estimated $\bL$, the diagonal elements can be
estimated by the sample variance of the components in
$\hat{L}\bw_i$. The approach can easily be adapted to estimate the
sparse precision matrix $\Sig^{-1}$.  See \cite{MB05}
for a similar approach and a thorough study.

\subsection{Health studies}

Survival data analysis has been a very active research topic because
survival data are frequently  collected from reliability analysis,
medical studies, and credit risks.
In practice, many covariates are often available as potential risk
factors.
Selecting significant variables plays a crucial role in model
building for survival data but is challenging due to the complicated
data structure.
Fan and Li \cite{FL02} derived the nonconcave penalized partial
likelihood for Cox's model and Cox's frailty model, the most
commonly used semiparametric models in survival analysis. Cai, \etal
\cite{CFLZ05} proposed a penalized pseudo partial likelihood for
marginal Cox's model with multivariate survival data and applied the
proposed methodology for a subset data in the Framingham study,
introduced in Section 2.2.

Let $T,\ C$ and $\bx$ be respectively the survival time, the
censoring time and their associated covariates. Correspondingly, let
$Z=\min\{T,C\}$ be the observed time and $\delta=I(T\le C)$ be the
censoring indicator. It is assumed that $T$ and $C$ are
conditionally independent given $\bx$,  that the censoring
mechanism is noninformative, and that the observed data $\{(\bx_i,
Z_i,\delta_i): i=1,\cdots,n\} $ is an independently and identically
distributed random sample from a certain population
$(\bx,Z,\delta)$.
%
%
The Cox  model assumes the conditional hazard function of $T$ given
$\bx$
\begin{equation}
h(t|\bx)=h_0(t)\exp(\bx^T\bbeta), \label{c1a}
\end{equation}
where   $h_0(t)$ is an unspecified
baseline hazard function.
Let $t_1^0<\cdots < t_N^0$ denote the ordered observed failure
times. Let $(j)$ provide the label for the item failing at $t_j^0$
so that the covariates associated with the $N$ failures are
$\bx_{(1)},\cdots, \bx_{(N)}$. Let $R_j$ denote the risk set right
before the time $t_{j}^0:$ $ R_j=\{i:Z_i\ge t_j^0\}. $
Fan and Li \cite{FL02} proposed the penalized partial likelihood
\beq Q(\bbeta)=\sum_{j=1}^N [ \bx_{(j)}^T\bbeta-\log\{\sum _{i\in
R_j} \exp (\bx_i^T\bbeta)\} ]
-n\sum_{j=1}^d \plj(|\beta_j|).
\label{c6a} \eeq
The penalized likelihood estimate of $\bbeta$ is to maximize
(\ref{c6a}) with respect to $\bbeta$.

For finite parameter settings, Fan and Li \cite{FL02} showed that
under certain regularity conditions, if both $a_n$ and $b_n$ tend to
0, then there exists a local maximizer $\hat{\bbeta}$ of the
penalized partial likelihood function in (\ref{c6a}) such that $\|
\hat{\bbeta} - \bbeta_0 \| = O_P(n^{-1/2} +a_n)$. They further
demonstrated the following oracle property.

\begin{theorem}  
Assume that  the penalty function $\pln(|\beta|)$ satisfies
condition (\ref{c3}). If $\lambda_{n,j}\to 0$,
$\sqrt{n}\lambda_{n,j}\to \infty$ and $a_n = O(n^{-1/2})$, then
under some mild regularity conditions, with probability tending to
1, the root $n$  consistent local maximizer $
\hat{\bbeta}=(\hat{\bbeta}_{1}^T, \hat{\bbeta}_{2}^T)^T $ of
$Q(\bbeta)$ defined in (\ref{c6a})  must satisfy
\[
\hat{\bbeta}_{2}={\bf 0},\quad\mbox{and}\quad
\sqrt{n}(I_1+\Sig)\left\{\hat{\bbeta}_{1}- \bbeta_{10} +
(I_1+\Sig)^{-1}\bb\right\} \wcon N\left\{{\bf
0},I_1(\bbeta_{10})\right\},
\]
where $I_1$ is the first $s \times s$ submatrix of the Fisher
information matrix $I(\bbeta_0)$ of the partial likelihood.
\end{theorem}

Cai, \etal \cite{CFLZ05} investigated the sampling properties of
penalized partial likelihood estimate with a diverging number of
predictors and clustered survival data. They showed that the oracle
property is still valid for penalized partial likelihood estimation
for the Cox marginal models with multivariate survival data.

\subsection{Financial engineering and risk management}

There are many outstanding challenges of dimensionality in diverse
fields of financial engineering and risk management.  To be concise,
we focus only on the issue of covariance matrix estimation using a
factor model.

Let $Y_i$ be the excess return of the $i$-th asset over the
risk-free asset.  Let $f_1, \cdots, f_K$ be the factors that
influence the returns of the market.  For example, in the
Fama-French 3-factor model, $f_1$, $f_2$ and $f_3$ are respectively
the excessive returns of the market porfolio, which is the
value-weighted return on all NYSE, AMEX and NASDAQ stocks over the
one-month Treasury bill rate, a portfolio constructed based on the
market capitalization, and a portfolio constructed based on the
book-to-market ratio.  Of course, constructing factors that
influence the market itself is a high-dimensional model selection
problem with massive amount of trading data.  The $K$-factor model
\cite{CRo83,Ros76} assumes
\begin{equation}
Y_i=b_{i1}f_1+\cdots+b_{iK}f_K+\varepsilon_i,\quad i=1,\cdots,d,
\label{f15}
\end{equation}
where $\{\varepsilon_i\}$ are idiosyncratic noises, uncorrelated
with the factors, and $d$ is the number of assets under
consideration. This an extension of the famous Capital Asset Pricing
Model derived by Sharpe and Lintner (See \cite{CLM97,Coc01}).
Putting it into the matrix form, we have
$\by=\bB\bff+\bveps$ so that
\begin{equation}
\Sig=\Var(\bB\bff) + \Var(\bveps) = \bB\Var(\bff)\bB^T +\Sig_0,
\label{f16}
\end{equation}
where $\Sig = \Var(\by)$ and $\Sig_0 = \Var(\bveps)$ is assumed to
be diagonal.

Suppose that we have observed the returns of $d$ stocks over $n$
periods (e.g., 3 years daily data).  Then, applying the
least-squares estimate separately to each stock in (\ref{f15}), we
obtain the estimates of coefficients in $\bB$ and $\Sig_0$.  Now,
estimating $\Var(\bff)$ by its sample variance, we obtain a
substitution estimator $\hSig$ using (\ref{f16}).  On the other
hand, we can also use the sample covariance matrix, denoted by
$\hSig_{\text{sam}}$, as an estimator.

In the risk management or portfolio allocation, the number of stocks
$d$ can be comparable with the sample size $n$ so it is better
modeled as $d_n$.  Fan \etal \cite{FFL05} investigated thoroughly
when the estimate $\hSig$ outperforms $\hSig_{\text{sam}}$ via both
asymptotic and simulation studies. Let us quote some of their
results.

\begin{theorem}   
Let $\lambda_k(\Sig)$ be the $k$-th largest eigenvalue of $\Sig$.
Then, under some regularity conditions, we have
$$
 \max_{1\leq k\leq d_n}\left|\lambda_k(\hSig)-\lambda_k(\Sig)\right|
    =  o_P \{ (\log n \; d_n^2/n)^{1/2} \} =
 \max_{1\leq k\leq
 d_n}\left|\lambda_k(\hSig_{\text{\rm sam}})-\lambda_k(\Sig)\right|.
$$
For a selected portfolio weight $\bxi_n$ with ${\bf 1}^T \bxi_n =1$,
we have
$$
\left|\bxi_n ^T \hSig \bxi_n - \bxi_n^T \Sig \bxi_n \right| =
o_P\{(\log n \; d_n^4 /n)^{1/2}\} =
   \left|\bxi_n^T \hSig_{\text{\rm sam}} \bxi_n-\bxi_n^T
   \Sig_n\bxi_n\right|.
$$
If, in addition, the all elements in $\bxi_n$ are positive, then the
latter rate can be replaced by $o_P\{( \log n \; d_n^2/n)^{1/2}\}$.
\end{theorem}

The above result shows that for risk management where the portfolio
risk is $ \bxi_n^T \Sig \bxi_n $, no substantial gain can be
realized by using the factor model.  Indeed, there is no substantial
gain for estimating the covariance matrix even if the factor model
is correct. These have also convincingly been demonstrated in
\cite{FFL05}
using simulation studies.  Fan \etal \cite{FFL05} also gives the
order $d_n$ under which the covariance matrix can be consistently
estimated.

The substantial gain can be realized if $\Sig^{-1}$ is estimated.
Hence, the factor model can be used to improve the construction of
the optimal mean-variance portfolio, which involves the inverse of
the covariance matrix.  Let us quote one theorem of \cite{FFL05}.
See other results therein for optimal portfolio allocation.

\begin{theorem}  
Under some regularity conditions, if $d_n = n^\alpha$, then for $0
\leq \alpha < 2$,
$$
  d_n^{-1} \mbox{\rm tr}(\Sig^{-1/2} \hSig
  \Sig^{-1/2} - I_{d_n})^2 = O_P(n^{-2\beta})
$$
with $\beta =  \min(1/2, 1-\alpha/2)$, whereas for $\alpha < 1$, $
d_n^{-1} \mbox{\rm tr}(\Sig^{-1/2} \hSig_{\text{sam}} \Sig^{-1/2} -
I_{d_n})^2 = O_P(d_n/n)$. In addition, under the Frobenius norm,
$$
d_n^2 \| \hSig^{-1} - \Sig^{-1} \|^2 = o(d_n^4 \log n/n) = \|
\hSig_{\text{sam}}^{-1} - \Sig^{-1} \|^2.
$$
\end{theorem}

\subsection{Machine learning and data mining}

In machine learning, our goal is to build a model with the
capability of good prediction of future observations. Prediction
error depends on the loss function, which is also referred to as a
divergence measure. Many loss functions are used in the literature.
To address the versatility of loss functions, let us
use the device introduced by \cite{Bre67}. 
For a concave function $q(\cdot)$, define a $q$-class of loss
function $\ell(\cdot,\cdot)$ to be
\begin{equation}
\ell(y,\hat{m}) = q(\hat{m}) -q(y) - q'(\hat{m})(\hat{m}-y)
\label{f12}
\end{equation}
where $\hat{m}\equiv \hat{m}(\bx)$, an estimate of the regression
function $ m(\bx)=E(y|\bx)$. Due to the concavity of $q$,
$\ell(\cdot,\cdot)$ is non-negative.

Here are some notable examples of $\ell$-loss constructed from the
$q$-function. For binary classification, $y\in \{-1,1\}$. Letting
$q(m) = 0.5\min\{1-m, 1+m\}$ yields the misclassification loss,
$\ell_1(y, \hat{m})=I\{y\ne  I (\hat{m}>0)\}$. Furthermore,
$\ell_2(y,\hat{m})= [1-y\sgn(\hat{m})]_+$ is the hinge loss if
$q(m)=\frac14 \min\{1-m, 1+m\}$. The function $q_3(m) =
\sqrt{1-m^2}$ results in $\ell_3(y,\hat{m})=\exp\{-0.5 y
\log\{(1+\hat{m})/(1-\hat{m})\}$, the exponential loss function in
AdaBoost \cite{FSc97}. Taking $q(m)=c m - m^2$ for some constant $c$
results in the quadratic loss $\ell_4(y,\hat{m}) = (y-\hat{m})^2$.

For a given loss function, we may extend the PLS to a penalized
empirical risk minimization (\ref{aaa}).
The dimensionality $d$ of the feature vectors can
be much larger than $n$ and hence the penalty is needed to select
important feature vectors.  See, for example, \cite{BL04}
for an important study in this direction.

We next make a connection between the penalized loss function and
the popularly used support vector machines (SVMs), which have been
successfully applied to various classification problems. In binary
classification problems, the response $y$ takes values either 1 or
$-1$, the class labels. A classification rule $\delta(\bx)$ is a
mapping from the feature vector $\bx$ to $\{1, -1\}$. Under the 0--1
loss, the misclassification error of $\delta$ is $P\{y\ne
\delta(\bx)\}.$ The smallest classification error is the Bayes error
achieved by $ \mbox{argmin}_{c\in \{1,-1\}} P(y=c|\bx). $ Let
$\{\bx_i,y_i\}$, $i=1,\cdots,n$ be a set of training data, where
$\bx_i$ is a vector with $d$ features, and the output $y_i\in
\{1,-1\}$ denotes the class label. The 2-norm SVM is to find a
hyperplane $\bx^T\bbeta $, in which $x_{i1}=1$ is an intercept and
$\bbeta=(\beta_1,\bbeta_{(2)}^T)^T$, that creates the biggest margin
between the training points from class 1 and $-1$
\cite{Va96}: \beq \max_{\bbeta} \frac 1{\|\bbeta_{(2)}\|^2}\qquad
 \mbox{subject\ to\ } y_i(\bbeta^T\bx_i)\ge 1-\xi_i, \forall\
i, \
 \xi_i \ge 0, \sum \xi_i \le B, \label{dd1}
\eeq where $\xi_i$ are slack variables, and $B$ is a pre-specified
positive number that controls the overlap between the two classes.
Due to its elegant margin interpretation and highly competitive
performance in practice, the 2-norm SVM has become popular and has
been applied for a number of classification problems. It is known
that the linear SVM has an equivalent hinge loss formulation
\cite{HTF01}
\[
\hbbeta = \mbox{argmin}_{\bbeta} \sum_{i=1}^n
[1-y_i(\bx_i^T\bbeta)]_+ + \lambda \sum_{j=2}^d \beta_j^2.
\]
Lin \cite{Lin02} shows that the SVM directly approximates the Bayes
rule without estimating the conditional class probability because of
the unique property of the hinge loss. As in the ridge regression,
the $L_2$-penalty helps control the model complexity to prevent
over-fitting.

Feature selection in the SVM has received increasing attention in
the literature of machine learning. For example, the last issue of
volume 3 (2002-2003) of {\em Journal of Machine Learning
Research} is a special issue on feature selection and extraction for
SVMs. We may consider a general penalized SVM
\[
\hbbeta = \mbox{argmin}_{\bbeta} \sum_{i=1}^n
[1-y_i(\bx_i^T\bbeta)]_+ + \sum_{j=1}^d \plj(|\beta_j|).
\]
%
%
The 1-norm (or LASSO-like) SVM has been used to accomplish the goal of automatic feature selection
in the SVM (\cite{ZRHT03}). 
Fried\-man \etal \cite{Fetal04} shows that the 1-norm SVM is
preferred if the underlying true model is sparse, while the 2-norm
SVM performs better if most of the predictors contribute to the
response.
With the SCAD penalty, the penalized SVM may improve the bias properties of the 1-norm SVM.



\def\aos{{\em Ann. Statist. }}
\def\bka{{\em Biometrika }}
\def\cjs{{\em Canad. J.  Statist. }}
\def\bcs{{\em Biometrics }}
\def\SS{{\em Statist.  Sinica }}
\def\bios{{\em Biostatistics }}
\def\aism{{\em Ann. Inst. Statist.  Math. }}
\def\jasa{{\em J. Amer.  Statist.  Assoc. }}
\def\sm{{\em Statist.  Med. }}
\def\jrssb{{\em J.  Royal Statist. Soc., B }}
\def\jrssc{{\em J.   Royal Statist.  Soc., C }}
\def\jcgs{{\em J. Comput. Graph.  Statist. }}
\def\jmva{{\em J.  Multiv.  Anal. }}


\end{document}